\documentclass[reqno, 12pt]{amsart}

\usepackage{amssymb}
\usepackage{latexsym}
\usepackage{epsfig}
\usepackage{epsf}
\usepackage{float}
\usepackage{a4wide}
\usepackage{hyperref}

\newtheorem{theo}{Theorem}
\newtheorem{lemm}{Lemma}

\newtheorem{rema}{Remark}

\def\R{{\mathbb R}}  
  
\def\P{{\mathbb P}}

\def\E{{\mathbb E}}

\def\0{{\mathbf 0}}

\newcommand{\Cov}{{\rm Cov}}

\def\proofof #1{{\noindent \emph{Proof of #1}.}}
\def\endproof{{\flushright$\square$ \vskip 2mm}}

\begin{document}

\title[On the location of the maximum of a continuous stochastic process]{On the location of the maximum of a continuous stochastic process}
\author{Leandro P. R. Pimentel}
\address{Instituto de Matem\'atica, Universidade Federal do Rio de Janeiro, Brasil.} 
\email{leandro@im.ufrj.br}



\begin{abstract}
In this short note we will provide a sufficient and necessary condition to have uniqueness of the location of the maximum of a stochastic process over an interval. The result will also express the mean value of the location in terms of the derivative of the expectation of the maximum of a linear perturbation of the underlying process. As an application, we will consider a Brownian motion with variable drift. The ideas behind the method of proof will also be useful to study the location of the maximum, over the real line, of a two-sided Brownian motion minus a parabola and of a stationary process minus a parabola.
\end{abstract}

\maketitle 

\section{Introduction}\label{sec:Introd}
Let $(X(z)\,,\,\,z\in [s,t])$ be a stochastic process with continuous paths on $[s,t] \subseteq\R$. The maximum of $X$ on $[s,t]$ is defined as 
$$M(X):=\max_{z\in [s,t]}X(z)\,,$$
 and the set of locations of the maximum (or $\arg\max$) is defined as
$$\arg\max(X):=\left\{z\in [s,t]\,:\,X(z)=M\right\}\,.$$
By continuity, $M$ is well defined and $\arg\max(X)$ is a nonempty compact subset of $[s,t]$. In many situations, we do expect that the maximum is actually attained almost surely at a unique location $Z$, so that 
\begin{equation}\label{uniqueness}
\arg\max(X)\stackrel{a.s.}{=}\{Z\}\,.
\end{equation}
In this article, we will prove a sufficient and necessary condition to have \eqref{uniqueness}. The main result is stated below.
\begin{theo}\label{theo:unique}
Let $(X(z)\,,\,\,z\in [s,t])$ be a stochastic process with continuous paths on $[s,t]$ and assume that  that $\E |M|<\infty$.  For $a\in\R$ let 
$$X^a(z):=X(z)+a z\,,\,\,z\in [s,t]\,,$$
and define
$$M^a:=M(X^a)\,,\,\mbox{ and }\,\,m(a):=\E M^a\,.$$
Then $a\mapsto m(a)$ is differentiable at $a=0$ if and only if the location of the maximum is almost surely unique \eqref{uniqueness}. In the latter case we have:
\begin{equation}\label{eq:uniqueformula} 
\E Z=m'(0)\,, 
\end{equation}
where $m'(0)$ is the derivative of $m$ at $a=0$. 
\end{theo}

The proof of Theorem \ref{theo:unique} is based on a simple non-probabilistic result (Lemma \ref{main}), which roughly states that the left and right directional derivatives of the functional $M$, with respect to the identity function, are given by the left-most  and the right-most locations of the maximum, respectively \footnote{Notice that if $X(z)=0$ for all $z\in[0,t]$ then $m(a)$ is not differentiable at $a=0$, since $m(a)=0$ for $a\leq0$ and $m(a)=a t$ for $a>0$.}.  

An example where we can apply Theorem \ref{theo:unique} is given by $X=B+f$, where $B$ is a Brownian motion and $f$ is a deterministic continuous function.    
\begin{theo}\label{cov}
Let 
$$X(z)=B(z)+f(z)\,,\,\,\,z\in[0,t]\,,$$
where $B$ is a standard Brownian motion process and $f$ is a deterministic continuous function. Then the location of the maximum is almost surely unique \eqref{uniqueness} and 
\begin{equation}\label{maxcov}
\E Z=\Cov\big(M,B(t)\big)\,. 
\end{equation}
\end{theo}      

Theorem \ref{cov} implies a similar result when $f$ is a continuous process that is independent of $B$. An interesting aspect of \eqref{maxcov} is that it gives the same result as if $Z$ were independent of $B$ \footnote{If $U\in[0,t]$ is independent of $B$ then $\Cov\big(B(U)+f(U),B(t)\big)=\Cov(B(U),B(t))= \E U$. Can we understand this behavior of $Z$, for $f\equiv 0$, in the light of Levy's $M-B$ theorem, or Pitman's $2M-B$ theorem?}.  Analogous identities have also appeared in particle systems and percolation models \cite{BaCaSe,CaGr,FeFo,Se}. The uniqueness of the location of the maximum for a continuous Gaussian process was proved by Kim and Pollard \cite{KiPo}, and it can certainly be used in our context. The author has tried to compute the derivative of $m$ for a Gausian process $X$, in order to provide an alternative proof of uniqueness based on Theorem \ref{theo:unique}, but with no success so far.       

In the previous situation we have considered the expectation of the maximum of a linear perturbation of the process $X$ and compute its derivative at zero. Other types of perturbation can also provide useful information about the location of the maximum. For instance, consider $X=B+f$, where now $B$ denotes a standard two-sided Brownian motion, $f$ is again a continuous function, and $z\in[-t,t]$. By taking a perturbation with respect to $z\mapsto z_+:=\max\{z,0\}$, we will see that
\begin{equation}\label{twosided}             
\E Z_+=\Cov\big(M,B(t)\big)\,.
\end{equation}

The maximum, over the real line, of a two-sided Brownian motion minus a parabola, and its location, arises as a limit object in many different statistical problems. Theorem \ref{theo:unique} can be used in this context as well to ensure uniqueness, since a.s. the $\arg\max$ will be compact (due to the negative parabolic drift). For many examples and various results see \cite{Gr}. By symmetry, is not hard to see that the location has zero mean. The expectation of the maximum and the variance of the location can be expressed in terms of integrals involving the Airy function \cite{Da,Gr1,Ja,Ja1}. Relying on those expressions, Groeneboom \cite{Gr1} and Janson \cite{Ja} remarked that the variance of the location equals one third the expectation of the maximum. By adding a quadratic perturbation (i.e. $a z^2$), and computing the derivative of the expected maximum, we can directly prove the Groeneboom-Janson relation. We also note that the probability that the maximum over the real line differs from the maximum over $[-t,t]$ decays exponentially fast to zero, and thus \eqref{twosided} can be extended to the limiting behavior.   
\begin{theo}\label{theo:JaGro}
Let 
$$X(z)=B(z)-z^2\,,\,\,\,z\in\R\,,$$
where $B$ is a standard two-sided Brownian motion process. Then the location of the maximum is almost surely unique \eqref{uniqueness} and 
\begin{equation}\label{JaGr}
\E Z=0\,,\,\,\E Z_+=\lim_{t\to\infty}\Cov\left(M,B(t)\right)\,\,\mbox{ and }\,\,\,\,\E Z^2 =\frac{1}{3}\E M\,.
\end{equation}
\end{theo}   

Another situation where uniqueness can be proved by the same methodology is when $X$ is a stationary process minus a parabola \footnote{Also when $X$ is Brownian motion minus a linear drift, over $[0,\infty)$. In this case, the distribution of  $M$ is well known to be exponential.}.
\begin{theo}\label{theo:stationary}
Let
$$X(z)=A(z)-z^2\,,\,\,\,z\in\R\,,$$
where $A$ is a stationary process with continuous paths. Assume that 
\begin{equation}\label{statassump}
\E |M|<\infty\,\,\mbox{ and }\,\,\int_0^{\infty}\P\left(\arg\max(X)\not\subseteq [-u,u]\right)du <\infty\,.
\end{equation}
Then the location of the maximum is almost surely unique \eqref{uniqueness} and 
\begin{equation}\label{statmean}
\E Z= 0\,.
\end{equation}
\end{theo}

It is surprising that \eqref{statmean} holds for any stationary process minus a parabola. However, as we shall see, the derivative of $m(a)$ can be easily computed in this case. The Airy process \cite{PrSp} is a an example where Theorem \ref{theo:stationary} can be used. It is a one-dimensional stationary process with continuous paths, whose finite dimensional distributions are described by Fredholm determinants. The interest in this process is mainly due to the fact that it gives the limit fluctuations of a number of processes appearing in statistical mechanics. Under the assumption that the maximum is indeed attained at a unique location, Johansson \cite{Jo} was able to prove that the law of the location describes the limit transversal fluctuations of maximal paths in last passage percolation models. This assumption was proved to be true by Corwin and Hammond \cite{CoHa}, and by Flores, Quastel and Remenik \cite{FlQuRe}. Both proofs used very strong results that depend on particular features of the Airy process. Theorem \ref{theo:stationary} is an alternative way to get uniqueness. 

\section{Proofs}\label{proof}

\subsection{Theorem \ref{theo:unique}}
Let $h:[s,t]\to\R$ be a continuous real function and let
$$Z_1(h):=\inf\arg\max(h)\,\,\mbox{ and }\,\,Z_2(h):=\sup\arg\max(h)\,.$$
We start with the analytic counterpart of the proof, that is given by Lemma \ref{main} below. It shows that the left and right directional derivatives of the functional $M$, with respect to the identity function, are given by $Z_1$ and $Z_2$, respectively.   
\begin{lemm}\label{main}
Let
$$h^a(z):=h(z)+a z\,.$$ 
Then
\begin{equation}\label{argmaxcont}
\lim_{a\to 0^-}Z_1(h^a)=Z_1(h)\,\,\mbox{ and }\,\,\lim_{a\to 0^+}Z_2(h^a)=Z_2(h)\,.
\end{equation}
Furthermore, 
\begin{equation}\label{maxder}
\lim_{a\to 0^-}\frac{M(h^a)-M(h)}{a}=Z_1(h)\,\,\mbox{ and }\,\,\lim_{a\to 0^+}\frac{M(h^a)-M(h)}{a}=Z_2(h)\,.
\end{equation}
\end{lemm}

\proof For simple notation, put $M^a=M(h^a)$ and $Z_i^a=Z_i(h^a)$. By continuity of $h$, 
\begin{equation}\label{maxbound}
M+a Z_i=h(Z_i)+a Z_i\leq M^a =h(Z^a_i)+a Z^a_i\leq M+a Z_i^a\,.
\end{equation}
This implies that 
\begin{equation}\label{maxderbound}
0\leq (M^a-M)-a Z_i \leq a(Z_i^a-Z_i)\,.
\end{equation}
The left-hand side inequality in \eqref{maxderbound} is equivalent to 
\begin{equation}\label{maxcont}
0\leq a(Z_i^a-Z_i)-(h(Z_i)-h(Z_i^a))\,.
\end{equation}
Since $h(Z_i)\geq h(Z_i^a)$, \eqref{maxcont} yields 
\begin{equation}\label{argmaxcont1}
Z_i^a\leq Z_i\,,\,\mbox{ for }\,\,a<0\,,\,\mbox{ and }\,\,Z_i^a\geq Z_i\,,\,\,\mbox{ for }\,\,a>0\,.
\end{equation}
By \eqref{argmaxcont1}, if \eqref{argmaxcont} is not true for $i=1$, then there exist $\delta>0$ and a sequence $a_n\to0^-$ such that  $Z_1^{a_n}\leq Z_1-\delta$ for all $n\geq 1$. By compactness of $[s,t]$, one can find a subsequence $a_{n_k}\to0^-$ and $\tilde Z_1\in K$ such that $\tilde Z_1=\lim_{k\to\infty}Z_1^{a_{n_k}}\leq Z_1-\delta$. By \eqref{maxcont} (and continuity of $h$), this implies that $h(\tilde Z_1)\geq h(Z_1)$, which leads to a contradiction, since $Z_1$ is the left-most location of the maximum. The proof for $i=2$ is analogous. 

Now, by \eqref{maxderbound},
\begin{equation}\label{eq:maxderiv1}
0\,\,\geq\,\, \frac{M^a-M}{a}- Z_1\,\,\geq \,\,Z_1^a-Z_1\geq s-t\,,\,\,\mbox{ for }\,\,a<0\,,
\end{equation}
and
\begin{equation}\label{eq:maxderiv2}
0\,\,\leq\,\, \frac{M^a-M}{a}- Z_2\,\,\leq \,\,Z_2^a-Z_2\leq t-s\,,\,\,\mbox{ for }\,\,a>0\,.
\end{equation}
Together with \eqref{argmaxcont}, \eqref{eq:maxderiv1}  and \eqref{eq:maxderiv2} imply  \eqref{maxder}. 
\endproof

\proofof{Theorem \ref{theo:unique}}
Notice that \eqref{maxbound} implies that $|M^a|$ has finite expectation, since we assume that  $\E|M|<\infty$, and $Z_i,Z_i^\epsilon\in [s,t]$. Also, the distance between $Z_i$ and $Z_i^a$ is always bounded by $t-s$. This will be important in the probabilistic counterpart of the proof, in order to use dominated convergence, as follows: If $m(a)$ is differentiable at $a=0$ then
$$m'(0)=\lim_{a\to 0^-}\frac{m(a)-m(0)}{a}=\lim_{a\to 0^+}\frac{m(a)-m(0)}{a}\,.$$
Together with \eqref{maxder}, and dominated convergence, this proves that $\E Z_1=\E Z_2$. Since $Z_1\leq Z_2$, we must have that $Z_1\stackrel{a.s.}{=}Z_2$, which yields to \eqref{uniqueness} and \eqref{eq:uniqueformula}. Reciprocally, if \eqref{uniqueness} is true then, by Lemma \ref{main}, 
$$\lim_{a\to 0^-}\frac{M^a-M}{a}\stackrel{a.s.}{=}\lim_{a\to 0^+}\frac{M^a-M}{a}\,.$$
Thus, dominated convergence implies that 
$$\lim_{a\to 0^-}\frac{m(a)-m(0)}{a}=\lim_{a\to 0^+}\frac{m(a)-m(0)}{a}\,,$$
which shows that $m(a)$ is differentiable at $a=0$, and the proof is finished.
\endproof

\subsection{Theorem \ref{cov}}
In the next lemma, we take $X=B+f$ and compute the derivative of $m(a)$ in a different way. This derivative can be computed by using the Cameron-Martin theorem. For sake of simplicity, we will present an alternative proof, which only requires basic knowledge of Brownian motion. 
\begin{lemm}\label{direcderiv} Let $Y=Y(B)$ be a (measurable) functional of standard Brownian motion $B$ on $[0,t]$ satisfying $\E Y^2<\infty$. Define  
$$y(a):=\E Y^a\,,\mbox{ where }\,\,Y^a:=Y(B^a)\,,$$
and $B^a(z):=az+B(z)$. Assume that $y(\cdot)$ is well defined in a neighborhood of $a=0$. Then
\begin{equation}\label{der}
y'(0)= \Cov(Y,B(t))\,.
\end{equation}
\end{lemm}

\proof Without loss of generality, we assume that $t=1$. The Brownian motion can be decomposed into
$$B(z)\stackrel{dist.}{=}N z+  B_0(z)$$
(as processes), where $B_0$ is a standard Brownian bridge with $B_0(0)=B_0(1)=0$, and $N$ is an independent Normal random variable of mean $0$ and variance $1$. Thus, 
$$B^a(z)=az+ B(z)\stackrel{dist.}{=}(a+N)z+ B_0(z)\,.$$
Since $B^a(1)=u$ if, and only if, $a+N=u$, we have that
\begin{equation}\label{bridge1}
B^a(z)\stackrel{dist.}{=}uz+ B_0(z)\,,
\end{equation}
conditioned on the event that $B^a(1)=u$. By \eqref{bridge1}, the conditional expectation of $Y^{a}$, given that $B^{a}(1)=u$, does not depend on $a\in\R$. Precisely, denote $B_u$ the process on the right hand side of \eqref{bridge1}. Then 
\begin{equation}\label{bridge2}
\E\left( Y^{a}\mid B^{a}(1)=u\right)=\E\left( Y\left(B_{u}\right)\right)\,.
\end{equation}
Therefore, by writing 
$$\rho_u(a):=\frac{1}{\sqrt{2\pi }}\exp\left\{-\frac{(u-a)^2}{2}\right\}\,,$$
we have that 
 $$y(a)=\int\E\left( Y^{a}\mid B^{a}(1)=u\right)\rho_u(a)du=\int\E\left( Y\left(B_{u}\right)\right)\rho_u(a)du\,.$$
Hence (by interchanging the derivative with the integral)
\begin{eqnarray}
\nonumber y'(a)&=&\int\E\left( Y^{a}\mid B^{a}(1)=u\right)\rho'_u(a)du\\
\nonumber&=&\int\E\left( Y^{a}\mid B^{a}(1)=u\right)\left(u-a \right)\rho_u(a)du\\
\nonumber&=& \E\left( Y^{a} B^{a}(1)\right)-a \E\left( Y^{a}\right)\,,
\end{eqnarray}
which proves \eqref{der}.
\endproof

\proofof{Theorem \ref{cov}} Take $Y(B):=M(B+f)$. By Lemma \ref{direcderiv},   
$$ m'(0)=\Cov(M,B(t))\,,$$
and hence, Theorem \ref{theo:unique} implies Theorem \ref{cov}. 
\endproof

\begin{rema}\label{remacov}
Given a square integrable function $\phi$ on $[0,t]$, define the function $\psi$ on $[0,t]$ by 
$$\psi(z):=\int_0^z\phi(u)du\,.$$
By the Cameron-Martin theorem, if $Y=Y(B)$ is a (measurable) functional of standard Brownian motion $B$ on $[0,t]$ satisfying $\E Y^2<\infty$ then  
$$\lim_{a\to 0}\frac{\E Y(B+a \psi)-\E Y}{a}=\E\left(Y\int_0^t\phi(z)dB(z)\right)\,.$$
If $\psi$ is increasing, then the same reasoning to prove Lemma \ref{main} yields to 
$$\lim_{a\to 0}\frac{M(h^{a,\psi})-M(h)}{a}=\psi(Z)\,$$
where $h^{a,\psi}(z):=h(a)+a\psi(z)$. Therefore
\begin{equation}\label{Cameron-Martin}
\E \psi(Z)=\E\left(M\int_0^t\phi(u)dB(u)\right)\,.
\end{equation}
By the chain rule, we also have that
\begin{equation}\label{chain}  
 \E\left(H'(M)Z\right)=\E\left(H(M)B(t)\right)\,.
\end{equation}
\end{rema}

\proofof{\eqref{twosided}} The proof is very similar. Put  $X^{a,+}(z):=X(z)+az_+$, for $z\in\R$, and $M^{a,+}=M(X^{a,+})$. Then 
$$\lim_{a\to 0^+}\frac{M^{a,+}-M}{a}=Z_{+}\,,$$
which implies that 
$$\lim_{a\to 0^+}\frac{m^+(a)-m^+(0)}{a}=\E Z_{+}\,,$$
where $m^+(a)=\E M^{a,+}$. By conditioning on $B(t)+a t$, this derivative equals  $\Cov(M,B(t))$, which shows \eqref{twosided}.
\endproof

\subsection{Theorem \ref{theo:JaGro}}
As we mentioned before, in this case, \eqref{uniqueness} can be obtained from Theorem \ref{cov} by using a.s. compactness of the $\arg\max$, and $\E Z=0$ follows easily from symmetry. The uniqueness also follows from Kim and Pollard \cite{KiPo}. For the sake of completeness, we also present an alternative proof which uses similar ideas as before. We start with a key lemma, which contains well known facts (see for instance Groeneboom \cite{Gr}):   
\begin{lemm}\label{lem2sided}
Let $B$ be a two sided Brownian motion and for $\beta\in\R$ define
$$M(\beta)=M(B,\beta):=\max_{z\in\R}\left\{B(z)-(z-\beta)^2 \right\}\,.$$ 
Let $Z_1(\beta)=Z_1(B,\beta)$ and $Z_2(\beta)=Z_2(B,\beta)$ denote the left-most and right-most locations of the maximum of $B(z)-(z-\beta)^2$, respectively. Then 
$$\E M(\beta)=\E M(0)\,\,\mbox{ and }\,\, \E Z_i(\beta)=\beta+\E Z_i(0)\,.$$
\end{lemm}

\proof By a.s. compactness of the $\arg\max$ (for $X(z)=B(z)-(z-\beta)^2$), $Z_i(\beta)$ is well defined for $i=1,2$. Notice that $\bar B\stackrel{dist.}{=}B$, where $\bar B(x):=B(x+\beta)-B(\beta)$ for $x\in\R$. Also,
$$\max_{z\in\R}\left\{B(z)-(z-\beta)^2 \right\}=\max_{x\in\R}\left\{\bar B(x)-x^2 \right\}+B(\beta)\,$$
(take $x=z-\beta$), and hence, 
$$M(B,\beta)=M(\bar B,0)+B(\beta)\,\,\mbox{ and }\,\,Z_i(B,\beta)-\beta=Z_i(\bar B,0)\,,$$
which proves the lemma. (Notice that the $\arg\max$ does not change by summing $B(\beta)$.)

\endproof 

\proofof{Theorem \ref{theo:JaGro}} By Lemma \ref{lem2sided}, 
$$m(a)=\E\max_{z\in\R}\left\{B(z)-z^2+a z \right\}=\E\max_{z\in\R}\left\{B(z)-\left(z-\frac{a}{2}\right)^2\right\}+\frac{a^2}{4}  = m(0)+\frac{a^2}{4}\,,$$
and hence 
\begin{equation}\label{der2sided}
m'(0)=0\,.
\end{equation}
Since the $\arg\max$ does not change by a vertical shifting of $a^2/4$, by Lemma \ref{lem2sided}, 
$$\E Z_i^a= \E Z_i(a/2)=a/2+\E Z_i\,,$$
which shows that 
\begin{equation}\label{distequal}
\lim_{a\to 0}\E Z_i^a = \E Z_i\,.
\end{equation}
On the other hand, by \eqref{maxderbound},
$$0\,\,\leq\,\, \frac{M^a-M}{a}- Z_i\,\,\leq \,\,Z_i^a-Z_i\,,\,\,\mbox{ for }\,\,a>0\,,$$
and
$$0\,\,\geq\,\, \frac{M^a-M}{a}- Z_i\,\,\geq \,\,Z_i^a-Z_i\,,\,\,\mbox{ for }\,\,a<0\,.$$
Together with \eqref{distequal}, these inequalities yield to 
$$\E Z_i=m'(0)\,,$$
and hence, $\E Z_1= \E Z_2$. Since $Z_1\leq Z_2$, we have that $Z_1\stackrel{a.s.}{=}Z_2$, which proves \eqref{uniqueness}. By \eqref{der2sided}, $\E Z=0$. To compute the limiting value of $\E Z_+$ use \eqref{twosided}.

To evaluate the second moment of $Z$, we add a quadratic perturbation to our original process and compute the derivative with respect to that. We follow the same notation as in \cite{Ja} and set $X_\gamma(z):=B(z)-\gamma z^2$ for $z\in\R$, 
$$M_\gamma:=M(X_\gamma)\,\,\mbox{ and }\,\,V_\gamma:=\arg\max(X_\gamma) \,\footnote{If the maximum is reached at a single value, then we refer to the point as the $\arg\max$.}\,.$$
Notice that 
$$M_{1-a}=\max_{z\in\R}\left\{B(z)- z^2+a z^2\right\}\,.$$
By scaling invariance of Brownian motion, 
\begin{equation}\label{scaling}
M_\gamma\stackrel{dist.}{=}\gamma_1^{1/3}\gamma^{-1/3}M_{\gamma_1}\,\,\mbox{ and }\,\,V_\gamma\stackrel{dist.}{=}\gamma_1^{2/3}\gamma^{-2/3}V_{\gamma_1}\,.
\end{equation}
Therefore,
$$n(a):=\E M_{1-a} =(1-a)^{-1/3} n(0)\,,$$
and thus,
\begin{equation}\label{dn1}
n'(0)=\frac{n(0)}{3}\,.
\end{equation}
On the other hand, as in the proof of \eqref{maxderbound}, 
$$0\leq (M_{1-a}-M_1)-a V_1^2 \leq a(V^2_{1-a}-V^2_1)\,,$$
which implies that, 
$$0\,\,\leq\,\, \frac{M_{1-a}-M_1}{a}- V^2_1\,\,\leq \,\,V^2_{1-a}-V^2_1\,,\,\,\mbox{ for }\,\,a>0\,,$$
and that
$$0\,\,\geq\,\, \frac{M_{1-a}-M_1}{a}- V^2_1\,\,\geq \,\,V^2_{1-a}-V^2_1\,,\,\,\mbox{ for }\,\,a<0\,.$$
By taking expectations on both sides of the last inequalities, and then using \eqref{scaling}, we have that   
$$\left |\frac{n(a)-n(0)}{a}-v(0)\right|\leq\left| v(a)-v(0)\right|=(1-a)^{-4/3}|v(0)|\,,$$
where $v(a)=\E V^2_{1-a}$.
Hence
$$n'(0)=v(0)\,.$$
Together with \eqref{dn1}, this shows that
$$\E Z^2=v(0)=n'(0)=\frac{n(0)}{3}=\frac{\E M}{3}\,.$$
We note that, by \eqref{scaling}, this also shows that $\E V^2_\gamma=(3\gamma)^{-1}\E M_\gamma$.
\endproof


\subsection{Theorem \ref{theo:stationary}} The proof of Theorem \ref{theo:stationary} is very similar to the previous one. 
\begin{lemm}\label{lem:stationary}
Let $A$ be a stationary process and for $\beta\in\R$ let 
$$M(\beta)=M(A,\beta):=\max_{s\in\R}\left\{A(s)-(s-\beta)^2 \right\}\,.$$ 
Let $Z_1(\beta)=Z_1(\beta)$ and $Z_2(\beta)=Z_2(\beta)$ denote the left-most and right-most locations of the maximum of $A(z)-(z-\beta)^2$, respectively. Then, for each fixed $\beta\in\R$, 
$$M(\beta)\stackrel{dist.}{=} M(0)\,\,\mbox{ and }\,\,  Z_i(\beta)-\beta\stackrel{dist.}{=} Z_i(0)\,.$$
\end{lemm}

\proof By stationarity, $\bar A\stackrel{dist.}{=}A$, where $\bar A(x):=A(x+\beta)$. On the other hand,
$$M(A,\beta)=\max_{z\in\R}\left\{A(z)-(z-\beta)^2 \right\}=M(\bar A,0)\,\,\mbox{ and }\,\,Z_i(A,\beta)-\beta=Z_i(\bar A,0)\,$$
(take $x=z-a$), which proves the lemma.
\endproof

\proofof{Theorem \ref{theo:stationary}} By Lemma \ref{lem:stationary},
$$m(a)=\E\max_{z\in\R}\left\{A(z)-z^2+a z \right\}=\E\max_{z\in\R}\left\{A(z)-\left(z-\frac{a}{2}\right)^2\right\}+\frac{a^2}{4}  = m(0)+\frac{a^2}{4}\,,$$
and hence,
\begin{equation}\label{eq:derstat}
m'(0)=0\,.
\end{equation}
By assumption \eqref{statassump}, we have that $\E |Z_i|<\infty$. Since the $\arg\max$ does not change by a vertical shifting of $a^2/4$, by Lemma \ref{lem:stationary}, 
$$\E Z_i^a= \E Z_i(a/2)=a/2+\E Z_i\,,$$
which shows that 
\begin{equation}\label{eq:distequal}
\lim_{a\to 0}\E Z_i^a= \E Z_i\,.
\end{equation}
On the other hand, by \eqref{maxderbound},
$$0\,\,\leq\,\, \frac{M^a-M}{a}- Z_i\,\,\leq \,\,Z_i^a-Z_i\,,\,\,\mbox{ for }\,\,a>0\,,$$
and
$$0\,\,\geq\,\, \frac{M^a-M}{a}- Z_i\,\,\geq \,\,Z_i^a-Z_i\,,\,\,\mbox{ for }\,\,a<0\,.$$
Together with \eqref{eq:distequal}, these inequalities yield to 
$$\E Z_i=m'(0)\,.$$
Thus, $\E Z_1= \E Z_2$. Since $Z_1\leq Z_2$, we have that $Z_1\stackrel{a.s.}{=}Z_2$, which proves \eqref{uniqueness}. By \eqref{eq:derstat}, $\E Z=0$.
\endproof

\paragraph{\bf Acknowledgements} Parts of this note were written during my visit to TU-Delft (December 2011) and my participation in the School and Workshop on Random Polymers and Related Topics at IMS-Singapore (May 2012), and I wish to thank the organizers of both events for their hospitality extended to me. During these periods I had the opportunity to meet many people and I would like to thank them for their collaboration, support and illuminating discussions. I also would like to say a special thanks to Eric Cator, Francesco Caravenna, and an anonymous referee, for the careful reading and helpful reviews of previous versions of this  manuscript.

\end{document}